\documentclass[12pt]{article}
\usepackage{amsfonts}

\newtheorem{Th.}{Theorem}[section]
\newtheorem{Lemma}[Th.]{Lemma}
\newtheorem{Proposition}[Th.]{Proposition}

\newtheorem{Corollary}[Th.]{Corollary}

\def\solidbox{\vrule width.6em height.5em depth.1em\relax}
\def\qed{\ifhmode\unskip\nobreak\fi\quad \solidbox}

\renewcommand{\theequation}{\arabic{section}.\arabic{equation}}
\newcommand{\msection}[1]{\section{#1}\setcounter{equation}{0}}

\newenvironment{meqnarray}[1]%
 {\\ \begin{minipage}{11cm}\begin{eqnarray*}\label{#1}}%
 {\end{eqnarray*}\end{minipage}\addtocounter{equation}{1}%
  \hfill (\theequation) \bigskip \\}

 {\begin{eqnarray*}\label{#1}}%
 {\end{eqnarray*}\addtocounter{equation}{1}}

\newenvironment{Proof}{\noindent{\bf Proof.}}%
                      {\hspace*{\fill} \solidbox  \bigskip \\}

\begin{document}


\title{The Balian--Low theorem for the symplectic form on ${\mathbb R}^{2d}$}
\author{John J. Benedetto,
        Wojciech Czaja, 
        Andrei Ya. Maltsev}
\date{\today}

\maketitle

\begin{abstract}
In this paper we extend the Balian--Low theorem, which is a version of the 
uncertainty principle for Gabor (Weyl--Heisenberg) systems, to functions of 
several variables.

In particular, we first prove the Balian--Low theorem for arbitrary quadratic 
forms. Then we generalize further and prove the Balian--Low theorem for 
differential operators associated with a symplectic basis for the symplectic 
form on ${\mathbb R}^{2d}$.

\end{abstract}

\newpage

\msection{Introduction}

For a given function $g \in L^2({\mathbb R}^d)$ we define the following two 
unitary operators on $L^2({\mathbb R}^d)$:
\[
M_{n} (g) (x) = e^{2\pi i n\cdot x} g(x), \quad n \in {\mathbb R}^d, 
\]
and
\[
T_{m} (g) (x) = g(x - m), \quad m \in {\mathbb R}^d,
\]
called {\it modulation} and {\it translation} operators, respectively.
In 1946 Dennis Gabor \cite{Ga} proposed to use these operators to define 
the collections of functions 
\[
g_{m,n} (x) = e^{2\pi i n\cdot x} g(x - m), 
\quad m,\,n \in {\mathbb Z},
\]
to be used in the analysis of information conveyed by communications channels.
These systems have been studied extensively in recent years. The edited
books by Benedetto and Frazier \cite{BF} and by 
Feichtinger and Strohmer \cite{FS}, as well as Gr\"ochenig's treatise 
\cite{G}, provide detailed treatments of various issues of the theory. Gabor systems 
are especially interesting because of their effective role in the time-frequency 
analysis of a wide variety of signals.

Let us now introduce some terms and notation that will be used throughout this 
paper. We say that a collection $\{ f_k : k=1, \ldots\}\subset 
L^2({\mathbb R}^d)$ of functions is a {\it frame} for $L^2({\mathbb R}^d)$, 
with {\it frame bounds} $A$ and $B$, if 
\[
\forall\; f \in L^2({\mathbb R}^d), \quad
A\| f \|_2^2 \le \sum_k |\langle f,f_k \rangle|^2 \le B\| f \|_2^2.
\]
A frame is {\it tight} if $A=B$; and a
frame is {\it exact} if it is no longer a frame after removal of any 
of its elements.
For any frame $\{ f_k : k=1, \ldots\}$ there exists a {\it dual frame}
$\{ {\tilde f}_k : k=1, \ldots\}$ such that
\begin{equation}
\label{e.*}
\forall\; f \in L^2({\mathbb R}^d), \quad
f = \sum_k \langle f,f_k \rangle {\tilde f}_k
= \sum_k \langle f, {\tilde f}_k \rangle f_k,
\end{equation}
where the series converge in $L^2({\mathbb R}^d)$. 
The choice of coefficients for expressing $f$ in terms of 
$\{ f_k : k=1, \ldots\}$ or $\{ {\tilde f}_k : k=1, \ldots\}$
is not unique, unless the frame is a {\it basis}.
{\it A frame is a basis if and only if it is exact}, e.g., \cite{BHW}.

For a frame $\{ f_k : k=1, \ldots\} \subset L^2({\mathbb R}^d)$ we define 
the associated {\it frame operator} $S$ on $L^2({\mathbb R}^d)$ by the rule,
\[
\forall\; f \in L^2({\mathbb R}^d), \quad
S(f) = \sum_k \langle f,f_k \rangle f_k.
\]
$S$ is a bounded and invertible map of $L^2({\mathbb R}^d)$ onto itself.
Given a frame $\{ f_k : k=1, \ldots\}$, our canonical choice of the dual 
frame $\{ {\tilde f}_k : k=1, \ldots\}$ will be defined by 
${\tilde f}_k  = S^{-1} (f_k)$. If a frame is exact then 
$\{ f_k : k=1, \ldots\}$ and $\{ {\tilde f}_k : k=1, \ldots\}$ are 
{\it biorthogonal}, that is,
\[
\langle f_k , {\tilde f}_l \rangle = \delta_{k,l} \quad k,l = 1, \ldots,
\]
where $\delta_{k,l}$ denotes the {\it Kronecker delta function}, i.e.,
it is 1 if $k=l$ and 0 otherwise.
It is elementary to show that $S^{-1}(g_{m,n}) = (S^{-1}(g))_{m,n}$
for Gabor frames $\{g_{m,n}\}$.

The {\it Fourier transform} is the unitary transformation $\cal F$
of $L^2({\mathbb R}^d)$ onto itself, defined formally by
\[
\hat{f} (\xi) = {\cal F} (f) (\xi) 
= \int_{{\mathbb R}^d} f(x) e^{2\pi i x\cdot \xi}
\;dx.
\]
We write ${\mathbb R}^d$ for arguments of a function $f\in L^2({\mathbb R}^d)$
and $\widehat{\mathbb R}^d$ for arguments of its Fourier transform.

We employ the standard notation in harmonic analysis, e.g., \cite{SW}.

The following result is a version
of the uncertainty principle for Gabor systems for the case $d=1$. 
It was first proved independently by Balian \cite{Bal} and Low \cite{L}. 
Both proofs contained a gap, which was corrected; 
and the result was generalized by Coifman, Daubechies, and Semmes 
from Gabor systems which form orthonormal bases to Gabor systems which form 
exact frames \cite{D}, see also \cite{B}, \cite{BHW}.
A different prof of Theorem \ref{blt1} was given by Battle \cite{Batt}.
Battle proved also an analogous result for wavelets \cite{Bat}.
 
\begin{Th.}{\bf Balian--Low theorem (BLT).}
\label{blt1}
Let $g\in L^2({\mathbb R})$ have the property that
$\{ g_{m,n}: \, m,n \in {\mathbb Z} \}$ is a Gabor orthonormal basis for
$L^2({\mathbb R})$. Then 
\begin{equation}
\label{e.1.2}
\left(\int_{\mathbb R} |g(x)|^2|x|^2 \;dx \right)
\left(\int_{\widehat{\mathbb R}} |\hat g(\xi)|^2|\xi|^2 \;d\xi\right) = \infty.
\end{equation}
\end{Th.}

\noindent {\bf Remark.}
Our original goal in this paper was to obtain a generalization of Theorem 
\ref{blt1} for functions of several variables.
In the process, and after having obtained some of our main results, 
we became aware of the work 
of Gr\"ochenig, Han, Heil, and Kutyniok \cite{GHHK},
in which the authors also extend the Balian--Low theorem to 
$d$-dimensions. Two of their fundamental results may be compared with our 
Theorem \ref{t1} and Theorem \ref{t6}. In fact, Theorem \ref{t6} is identical 
with the BLT for non-lattices in \cite{GHHK} and Theorem \ref{t1} extends the 
weak BLT for lattices in \cite{GHHK} to more general position and momentum 
operators.
Further, using techniques from the theory of metaplectic representations,
the authors in \cite{GHHK} generalize Theorem \ref{t6} to a Balian--Low type 
theorem for exact frames on symplectic lattices; for their setting their 
assertion states that there exists $i\in \{1, \ldots, d\}$ such that (\ref{e0}) 
below holds.

We follow a different path and prove that the choice of coordinates 
in (\ref{e0}) is not canonical, i.e., there is no ``preference"
for the directional derivatives and for multiplications by the standard 
basis coordinates. This means that one can work in any representation of 
${\mathbb R}^d$, e.g., Theorem \ref{t12}.
\hskip0.6cm

In Section 2 we prove the generalization 
of the Balian--Low theorem to $d$-dimensions in the standard coordinate system;
this is Theorem \ref{t1}.
As a corollary, we prove a Balian--Low theorem for arbitrary non-negative
quadratic forms (Corollary \ref{t3}).
In Section 3 we state and prove our main results, Theorem \ref{t12} and 
Theorem \ref{t13}, 
which assert a Balian--Low phenomenon (\ref{e.10*}) similar to but more 
far-reaching than (\ref{e.1.2}). The proof depends on our definition of 
generalized Fourier transforms which, in turn, allows us to reduce a rather 
general and comprehensive problem to the Balian--Low theorem in the standard 
coordinates as formulated in Theorem \ref{t6}.

Our approach is both straightforward and natural.
This is an essential part of our contribution.
It is also based on the quantum mechanical point of view.

\msection{Balian--Low theorem in standard coordinates}

Let ${ v}, { w} \in {\mathbb R}^d$ be non-zero vectors. We define 
the following operators, wherever they make sense in 
$L^2({\mathbb R}^d)$:
\[
P_{v} (f)(x) = \left(\sum_{i=1}^d {v}^i x^i\right) f(x) 
\]
and
\[
M_{w} (f)(x) 
= {\cal F}^{-1}\left(\left(\sum_{i=1}^d {w}^i \xi^i\right) \hat f (\xi)\right)(x) 
= {\cal F}^{-1}(P_{w} (\hat f)) (x),
\]
where $v = (v^1, \ldots, v^d) = \sum v^j u_j$, $u_j = (0,\ldots,0,1,0,\ldots,0)$ 
with $1$ in the $j$th coordinate, and $v^j \in {\mathbb R}$.
These unit vectors $u_j$ define the {\it standard Euclidean basis} 
$\{u_j:j=1,\ldots, d\}$ of ${\mathbb R}^d$.
If the vectors $v$ and $w$ in the definitions of $P_v$ and $M_w$ 
are elements of the standard basis, then we shall use 
the notation $P_i$ and $M_i$ for the operators induced by the {\it i}th basis
vector $u_i$.

The following result is our first generalization of the Balian--Low theorem.
The technique of proof is a well-known method for  
proving Balian--Low type theorems.
$\tilde g$ denotes the canonical dual defined in Section 1.

\begin{Th.}
\label{t1}
Let $\{ g_{m,n} : m,n \in {\mathbb Z}^d \}$ be an exact frame for 
$L^2({\mathbb R}^d)$. If $v, \, w \in {\mathbb R}^d$ 
satisfy $v \cdot w \neq 0$, then
\begin{equation}
\label{e.**}
\left\| P_{v}(g) \right\|_2 \left\| 
M_{w}(g) \right\|_2 \left\|P_{v} ({\tilde g}) \right\|_2 
\| M_{w}(\tilde g) \|_2 = \infty.
\end{equation}
\end{Th.}

\begin{Proof} 
We may assume without loss of generality that $|v| = |w| =1$, where 
$|\;|$ denotes the Euclidean norm in ${\mathbb R}^d$.
We shall proceed with a proof by contradiction; and so we assume 
that all four functions in (\ref{e.**}) are elements of $L^2({\mathbb R}^d)$.
Because of the biorthogonality relations for $g$ and $\tilde g$ we compute
\begin{meqnarray}{e.2.*}
\langle P_{v} (g), {\tilde g}_{{ m},{ n}} \rangle 
&=& \langle P_{v} (g), {\tilde g}_{m,n} \rangle 
- \left(\sum_{i=1}^d {v^i} m^i\right) \langle g, {\tilde g}_{m,n} \rangle \\
&=& \int_{{\mathbb R}^d} \left(\sum_{i=1}^d {v}^i (x^i - m^i)\right) 
g(x) \overline{{\tilde g}(x-m)} e^{-2\pi i n\cdot x}\; dx\\
&=& e^{-2\pi i m\cdot n} \int_{{\mathbb R}^d} \left(\sum_{i=1}^d {v}^i x^i\right)
\overline{{\tilde g}(x)} g(x+m)e^{-2\pi i n\cdot x}\; dx\\
&=& e^{-2\pi i m\cdot n} \langle g_{-m,-n}, P_{v}({\tilde g}) \rangle .
\end{meqnarray} 

From our assumption that $M_{ w} (g) \in L^2({\mathbb R}^d)$, it follows that 
the distributional partial derivative of $g$, $\partial_{w} g$, belongs to 
$L^2({\mathbb R}^d)$. From a standard result about Sobolev spaces, 
see, e.g., \cite{M}, Theorem 1.1,
there exists a function $h$ such that $g = h$ $a.e.$, and
$h$ is absolutely continuous on almost all straight lines parallel to the 
vector ${w}$. Thus the distributional directional derivative of $g$
coincides with the classical directional derivative $D_{w} (g)$ $a.e.$, and so
\[
M_{w} (g) (x) = \frac{i}{2\pi} D_{w} (g) (x) \quad a.e.
\]
Moreover, our assumptions imply that 
$D_{w} (g)$, $D_{w} (\tilde g) \in L^2({\mathbb R}^d)$. 
Therefore, using integration by parts, an appropriate change of variables, and 
the biorthogonality relations between $g$ and $\tilde g$, we can compute
\begin{meqnarray}{e.2.**}
\langle g_{m,n}, M_{w}({\tilde g}) \rangle 
\!\!&=& \frac{1}{2\pi i} \int_{{\mathbb R}^d} g(x-m) e^{2\pi i n\cdot x} 
\overline{D_{w} ({\tilde g})(x)}\;dx \\
&=& \frac{i}{2\pi} \int_{{\mathbb R}^d} D_{w} \left(g(x-m) 
e^{2\pi i n\cdot x}\right) \overline{{\tilde g}(x)} \; dx\\
&=& \frac{i}{2\pi} \int_{{\mathbb R}^d} \left( D_{w} (g)(x-m) 
e^{2\pi i n\cdot x} \right. \\
&+& \left. ({w} \cdot {n}) g(x-m) e^{2\pi i n\cdot x} \right) 
\overline{{\tilde g}(x)} \; dx\\
&=& \frac{i e^{2\pi i m\cdot n}}{2\pi} \! \int_{{\mathbb R}^d} \!\!
\left( D_{w} (g)(x) + 2\pi i ({w} \cdot {n}) g(x) \right)
e^{2\pi i n\cdot x}\overline{{\tilde g}(x+m)}  \; dx\\
&=& e^{2\pi i m\cdot n} \left( \langle  M_{w}(g), {\tilde g}_{-m,-n}\rangle 
+ ({w} \cdot n)\delta_{m,0} \delta_{n,0} \right)\\
&=& e^{2\pi i m\cdot n}\langle  M_{w}(g), {\tilde g}_{-m,-n}\rangle .
\end{meqnarray}

Because of (\ref{e.2.*}), (\ref{e.2.**}), and the frame representation property 
(\ref{e.*}), 
we have
\begin{meqnarray}{e.3*}
\langle P_v(g), M_w(\tilde g) \rangle 
&=& \sum_{m,n \in {\mathbb Z}^d} \langle P_v (g), {\tilde g}_{m,n} \rangle 
\langle g_{m,n}, M_w({\tilde g}) \rangle \\
&=& \sum_{m,n\in {\mathbb Z}^d} \langle  g_{-m,-n}, P_v({\tilde g}) \rangle
\langle  M_w(g), {\tilde g}_{-m,-n}\rangle \\
&=& \sum_{m,n\in {\mathbb Z}^d} \langle  M_w(g), {\tilde g}_{m,n}\rangle
\langle  g_{m,n}, P_v({\tilde g}) \rangle\\
&=& \langle M_w(g), P_v({\tilde g}) \rangle. 
\end{meqnarray}

It is not difficult to verify that
\begin{equation}
\label{e.2.***}
[P_{v}, M_{w}] = \frac{1}{2\pi i} ({v}\cdot {w}) \;{\rm Id},
\end{equation}
where the commutator $[P_{v}, M_{w}] = P_{v} M_{w} - M_{w}P_{v}$ and where
${\rm Id}$ denotes the identity operator, e.g., \cite{Mes} where 
(\ref{e.2.***}) appears for the position and momentum operators associated with
the standard basis vectors; see also the trivial calculation in \cite{BHW}. 
Thus, for functions $g$, $\tilde g \in L^2({\mathbb R}^d)$, such that
$P_{ v}(g)$, $P_{ v}(\tilde g) \in L^2({\mathbb R}^d)$ and 
$M_{ w} (g)$, $M_{ w}(\tilde g) \in L^2({\mathbb R}^d)$, 
we have
\begin{eqnarray*}
\langle P_{v}(g), M_{w}(\tilde g) \rangle &=& 
\langle M_{w}(g), P_{v}(\tilde g) \rangle
+ \frac{1}{2\pi i} ({v}\cdot {w}) \langle g, \tilde g \rangle\\
&=& \langle M_{ w}(g), P_{v}(\tilde g) \rangle
+ \frac{1}{2\pi i} ({v}\cdot {w}).
\end{eqnarray*}

Since we have assumed that ${ v}\cdot { w} \neq 0$, we obtain 
a contradiction with our calculation (\ref{e.3*}). 
\end{Proof}

\noindent{\bf Remark.} The claim (\ref{e.**}) is true if in Theorem \ref{t1} 
we consider the more general system $\{ g_{m,n} : (m,n) \in \Lambda\}$, where 
$\Lambda$ is an arbitraty lattice in ${\mathbb R}^{2d}$. For an analogous result
for position and momentum operators associated with the integer lattice 
${\mathbb Z}^{2d}$ see Theorem 8 in \cite{GHHK}.

\begin{Corollary}
\label{t2}
Let $\{ g_{m,n} : m,n \in {\mathbb Z}^d\}$ be an exact frame for 
$L^2({\mathbb R}^d)$. If ${v}, \, {w} \in {\mathbb R}^d$ 
satisfy ${v} \cdot {w} \neq 0$, then
\[
\left\| P_{v}(g) \right\|_2 \left\| 
M_{w}(g) \right\|_2 = \infty.
\]
\end{Corollary}

\begin{Proof} 
In view of Theorem \ref{t1}, it is enough to show that 
$P_{v}(g) \in L^2({\mathbb R}^d)$ if and only if
$P_{v}(\tilde g) \in L^2({\mathbb R}^d)$, and that 
$M_{w}(g)\in L^2({\mathbb R}^d)$
if and only if $M_{w}(\tilde g)\in L^2({\mathbb R}^d)$. 
This, in turn, was proved by Daubechies and Janssen \cite{DJ} for the position 
and momentum operators associated with the standard basis vectors, 
see also \cite{BHW}, Theorem 7.7.
The proof for arbitrary operators $P_v$ and $M_w$ is analogous, and it uses 
the $d$-dimensional Sobolev space argument which we have used in the proof 
of Theorem \ref{t1} instead of 1-dimensional considerations .
\end{Proof}

\noindent{\bf Example.} To show that the condition $v \cdot w \neq 0$ is 
necessary consider $L^2({\mathbb R}^2)$ with the orthonormal Gabor basis 
generated by 
\[
g(x,y) = \chi_{[0,1]} (x) \; {\cal F}^{-1}\left(\chi_{[0,1]}\right)(y)
\] 
and the vectors $v = (1,0)$ and $w = (0,1)$.
Then
\begin{eqnarray*}
\| P_v(g)\|_2^2 &=& \int_{{\mathbb R}^2} |x g(x,y)|^2\; dxdy \\
&=& \int_{\mathbb R} \left|x\chi_{[0,1]} (x)\right|^2\;dx 
\int_{\mathbb R}\left|{\cal F}^{-1}(\chi_{[0,1]})(y)\right|^2\;dy < \infty
\end{eqnarray*}
and
\begin{eqnarray*}
\| M_w(g)\|_2^2 &=& \int_{{\mathbb R}^2} |\eta \hat g(\xi,\eta)|^2\;d\xi d\eta\\
&=& \int_{\mathbb R}\left|{\cal F}(\chi_{[0,1]})(\xi)\right|^2\; d\xi
\int_{\mathbb R}\left|\eta (\chi_{[0,1]})(\eta)\right|^2\;d\eta< \infty.
\end{eqnarray*}
\vskip0.6cm

\begin{Corollary}
\label{t3}
Let $\omega (x)$ be any positive quadratic form on ${\mathbb R}^d$
and let $\{ g_{m,n} : m,n 
\in {\mathbb Z}^d\}$ be an exact frame for 
$L^2({\mathbb R}^d)$. Then
\[
\left( \int_{{\mathbb R}^d} \omega(x) |g(x)|^2\; dx \right) 
\left( \int_{{\widehat{\mathbb R}}^d} \omega(\xi) |\hat g(\xi)|^2\; d\xi \right)
= \infty.
\]
\end{Corollary}

\begin{Proof} 
Clearly, for any vector ${v}\neq 0$ we have 
${v} \cdot {v} \neq 0$. Thus, from 
Corollary \ref{t2} it follows that for any $\alpha_k \ge 0$ and 
$v_k \in {\mathbb R}^d$, $k=1, \ldots, d$, 
where some $\alpha_k > 0$, either
\[
\left(\sum_{k=1}^d \alpha_k \int_{{\mathbb R}^d} 
\left(\sum_{i=1}^d {v}^i_k x^i\right)^2 
|g(x)|^2\; dx \right) = \infty
\]
or
\[
\left( \sum_{k=1}^d \alpha_k \int_{{\widehat{\mathbb R}}^d}
\left(\sum_{i=1}^d {v}^i_k \xi^i\right)^2 
|\hat g(\xi)|^2\; d\xi \right) = \infty.
\]
The result follows since any quadratic form on ${\mathbb R}^d$ is of the form
\[
\omega(x) = \sum_{k=1}^d \alpha_k \left(\sum_{i=1}^d {v}^i_k x^i\right)^2,
\] 
where the $\alpha_k$s are non-negative.
\end{Proof}

We now consider a countable collection 
$\Lambda$ of vectors in ${\mathbb R}^{2d}$. 
For any pair $(m,n) \in \Lambda$, $m,n \in {\mathbb R}^d$, we shall associate 
the {\it translation-modulation transformation} $T_{m,n}$ defined on 
$L^2({\mathbb R}^d)$ as follows:
\[
T_{m,n} (g) (x) = e^{2 \pi i n \cdot x}
g\left(x+ m\right).
\]
From now on we shall write $g_{m,n} = T_{m,n} (g)$.
The study of {\it nonuniform} Gabor systems, i.e., those Gabor systems which 
are associated with a set $\Lambda$ which is not a lattice, has increased 
in recent years because of applications of such systems
to problems in signal processing, e.g., \cite{BPW}, \cite{BT}, 
\cite{Gro}, \cite{Li}. Of course, not all $\Lambda$s
generate orthonormal bases or even frames. In order for a Gabor system 
to have good signal representation properties, $\Lambda$ must satisfy certain 
density conditions. The most general results so far in this direction were 
obtained by Ramanathan and Steger \cite{RS} and by Christensen, Deng, 
and Heil \cite{CDH}.  

\vskip0.6cm
\noindent {\bf Example.} One easily constructs examples of 
uniform orthonormal Gabor bases for $L^2({\mathbb R}^d)$. 
The most simple example is $g(x) = \chi_{[0,1]^d}$ with the lattice 
$\Lambda = {\mathbb Z}^{2d}$. More interestingly there is the work
of Liu and Wang \cite{LW}, where the authors provide examples of 
nonuniform Gabor bases and frames, i.e., examples where $\Lambda$ 
is not a lattice.

For $d=1$ let $\Omega = [0,1] \cup [3,4]$ and 
\[
\Lambda = \left\{ 6{\mathbb Z}
+ \{-1,0,1\} \right\} \times \left\{ \frac{1}{2}{\mathbb Z} \right\}.
\]
Then $g(x) = (1/\sqrt{2}) \chi_{\Omega} (x)$ forms an orthonormal basis with 
translations and modulations in $\Lambda$. We would like to stress that 
although $\Lambda$ is a periodic set it is not a lattice, since in general 
a sum of two vectors in $\Lambda$ is not an element of $\Lambda$. 
We note that $\Lambda = -\Lambda$.

\cite{LW} also provides an account of 
differences between nonuniform Gabor bases in 1 and higher dimensions.
\vskip0.6cm

To prove Theorem \ref{t6} we shall need the following lemma, 
the proof of which is similar to the proof of analogous statements in Theorem 
\ref{t1}. 

\begin{Lemma}
\label{t7}
Let $g \in L^2({\mathbb R}^d)$, and let 
$\{{g}_{m,n} : (m,n) \in \Lambda \}$ be an orthonormal basis for 
$L^2({\mathbb R}^d)$. If $P_i (g), \, M_i (g) \in L^2({\mathbb R}^d)$, 
then 
\begin{equation}
\label{e1}
\langle {g}_{m,n}, P_i (g) \rangle 
= e^{2\pi i m\cdot n} \langle P_i (g), {g}_{-m,-n} \rangle  
\end{equation}
and
\begin{equation}
\label{e2}
\langle {g}_{m,n}, M_i (g) \rangle 
= e^{2\pi i m \cdot n} \langle M_i (g) , {g}_{-m,-n} \rangle.  
\end{equation}
\end{Lemma}

\begin{Proof} 
Since $\Lambda$ does not posses a lattice structure we cannot use 
(\ref{e.2.*}) and (\ref{e.2.**}). Indeed, the fact that a dual to a Gabor 
frame is also a frame of Gabor type holds only for systems associated with 
lattices. However, the assumption that
$\{ g_{m,n} : (m,n) \in \Lambda \}$ is an orthonormal basis for 
$L^2({\mathbb R}^d)$ compensates for this lack of structure in $\Lambda$.

\begin{eqnarray*}
\langle g_{m,n}, P_i (g) \rangle &=& \int_{{\mathbb R}^d} g(x-m) 
e^{2\pi i n\cdot x} \overline{P_i(g)(x)} \; dx \\
&=& e^{2\pi i m\cdot n} \int_{{\mathbb R}^d} g(x) 
e^{2\pi i n\cdot x} \overline{(x_i+m_i)g(x+m)} \; dx\\
&=& e^{2\pi i m\cdot n}\left( \langle P_i (g), {g}_{-m,-n} \rangle
+ m_i \langle g, {g}_{-m,-n} \rangle\right) \\
&=& e^{2\pi i m\cdot n} \langle P_i (g), {g}_{-m,-n} \rangle.
\end{eqnarray*}
The last equality above follows from the orthogonality of 
$\{ g_{m,n} : (m,n) \in \Lambda \}$. Similarly, using orthogonality and 
the integration by parts formula, we calculate

\begin{eqnarray*}
\langle {g}_{m,n}, M_i (g) \rangle &=& -\frac{i}{2\pi} 
\int_{{\mathbb R}^d} g(x-m)e^{2\pi i n\cdot x} 
\overline{D_i(g)(x)} \; dx \\
&=&\frac{i}{2\pi} \int_{{\mathbb R}^d}D_i\left( g(x-m) e^{2\pi i n\cdot x} 
\right) \overline{g(x)}\; dx\\
&=& \frac{i e^{2\pi i m \cdot n}}{2\pi} \int_{{\mathbb R}^d} \left( 
D_i(g)(x) + 2\pi im_i g(x) \right) e^{2\pi i n\cdot x} \overline{g(x+m)}\; dx \\
&=& e^{2\pi i m \cdot n} \left( \langle M_i (g) , g_{-m,-n} \rangle
+  m_i \langle g, g_{-m,-n} \rangle \right) \\
&=& e^{2\pi i m \cdot n} \langle M_i (g) , g_{-m,-n} \rangle .
\end{eqnarray*}

\end{Proof}

\begin{Th.}
\label{t6}
Let $\Lambda \subset{\mathbb R}^{2d}$ be a countable sequence of vectors 
with the property that $\Lambda = -\Lambda$.
Let $\{ T_{m,n}: (m,n) \in \Lambda\}$ be the 
associated family of translation-modulation transformations, and assume
$\{g_{m,n} :(m,n) \in \Lambda\}$ is an orthonormal basis for 
$L^2({\mathbb R}^d)$ for some $g\in L^2({\mathbb R}^d)$. 
For any $i = 1, \ldots, d$,
\begin{equation}
\label{e0}
\| P_i (g) \|_2 \| M_i (g) \|_2 = \infty .
\end{equation}
\end{Th.}

\begin{Proof} 
Because of (\ref{e1}), (\ref{e2}), the representation property of bases,
and the fact that $\Lambda = - \Lambda$, we obtain
\begin{meqnarray}{e.4*}
\langle M_i (g), P_i (g) \rangle &=& \sum_{(m,n) \in \Lambda}
\langle M_i (g) , g_{m,n} \rangle 
\langle g_{m,n}, P_i(g) \rangle\\
&=& \sum_{(m,n) \in \Lambda} \langle g_{-m,-n}, M_i (g) \rangle 
\langle P_i (g), g_{-m,-n} \rangle \\
&=& \langle P_i (g), M_i (g)\rangle .
\end{meqnarray}

On the other hand, again using the classical result from \cite{M}
used in Theorem \ref{t1}, we note that 
$M_i (g) \in L^2({\mathbb R}^d)$ implies that $\partial g/\partial x^i$ exists 
$a.e.$ Thus, integration by parts yields
\[
\langle M_i (g), P_i (g) \rangle 
= \langle P_i (g), M_i (g)\rangle - \frac{1}{2\pi i},
\]
which, in turn, leads to a contradiction with the calculation (\ref{e.4*}). 
\end{Proof}

\msection{Balian--Low theorem and symplectic forms}

The standard {\it symplectic form} $\Omega$ on ${\mathbb R}^{2d}$ is
defined as
\[
\Omega((x,y), (\xi, \eta)) = x\cdot \eta - y \cdot \xi,
\] 
for any $x,\, y,\, \xi, \, \eta \in {\mathbb R}^d$.
Note that $\Omega((x,0), (0,\xi)) = x \cdot \xi$. This observation, 
when compared to Theorem \ref{t1}, suggests a direction which we are going 
to follow in this section, and which yields our main result, Theorem \ref{t12}.

\vskip0.6cm
\noindent {\bf Definition.}
{\bf a.} A {\it symplectic basis} for ${\mathbb R}^{2d}$
with respect to the symplectic form $\Omega$ is a basis 
$\{a_j, b_j:j=1,\ldots, d\}\subset {\mathbb R}^{2d}$ for ${\mathbb R}^{2d}$
for which
\[
\Omega(a_i, a_j)= \Omega (b_i, b_j) =0
\]
and
\[
\Omega(a_i, b_j)= \delta_{i,j},
\]
for all $i, j =1, \ldots, d$. 

{\bf b.} If $\{\nu_i:i=1,\ldots,d\}\subset {\mathbb R}^d$ is any orthonormal 
basis for ${\mathbb R}^d$, then $a_i = (\nu_i,0)$, $b_i = (0,\nu_i)$, 
$i=1,\ldots,d$, is a symplectic basis for ${\mathbb R}^{2d}$.
For a non-trivial example in ${\mathbb R}^4$ take the row vectors of the matrix
\[
\left(\begin{array}{cccc} 
1 & 0 & - \frac{\sqrt{3}}{2} & - \frac{1}{2}\\
\frac{\sqrt{2}}{2} & - \frac{\sqrt{6}}{2} & 0 & - \frac{\sqrt{2}}{2}\\
0 & 1 & \frac{1}{2} & \frac{\sqrt{3}}{2} \\
- \frac{\sqrt{6}}{2} & \frac{\sqrt{2}}{2} & \frac{\sqrt{2}}{2} & 0
\end{array}\right).
\]

{\bf c.}
Consider the space ${\mathbb R}^{2d}$
with coordinates $(x^1, \ldots, x^d, y^1, \ldots, y^d)$
and let $\Omega$ be the symplectic form on ${\mathbb R}^{2d}$.
A {\it Lagrangian plane} $\Pi$ in ${\mathbb R}^{2d}$ is a $d$-dimensional 
subspace with the property that 
\[
\Omega|_{\Pi} =0.
\] 
\vskip0.6cm

If $\Pi$ is a Lagrangian plane in ${\mathbb R}^{2d}$ 
and if ${v}_1, \ldots, {v}_d \in \Pi \subset {\mathbb R}^{2d}$ is a basis
for $\Pi$ then, in particular, we have
\[
\Omega({v}_i, {v}_j) =0,
\]
for all $i,j = 1, \ldots, d$. For classical treatments of these and other 
related notions see, e.g., \cite{A}, \cite{AN}. 
A similar approach is used by H\"ormander \cite{Hor} to define Fourier Integral 
Operators, a special case of which we consider below.
A recent exposition of related results in case of Hermitian symplectic geometry
is due to Harmer \cite{H}.

We now define the differential operators 
$\{Q_{{v}_j}, \, j = 1, \ldots, d\}$ associated with a given basis
$\{v_j: j=1,\ldots,d\}$ for a given Lagrangian plane $\Pi$. 
Each $Q_{{v}_j}$ is defined by its action on a function $h$ as follows:
\begin{equation}
\label{eQ}
Q_{{ v}_j} (h)(x)= \frac{i}{2\pi} \nabla_j (h)(x) + f_j(x) h(x),
\end{equation}
where 
\[
\nabla_j = \sum_{k=1}^d v^{k+d}_j \frac{\partial}{\partial x^k} 
\]
and
\[
f_j(x)= \sum_{k=1}^d v^k_j x^k.
\]
Recall that $v_j^k$ is the $k$th coordinate of the vector 
$v_j \in {\mathbb R}^{2d}$ and that 
$x = \sum_{k=1}^d x^k u_k \in {\mathbb R}^d$.

The next result serves as the main motivation for our
work. It is analogous to a similar observation about commutators of
position and momentum operators that was asserted in equation (\ref{e.2.***}).
Its proof is also a straightforward calculation.

\begin{Proposition}
\label{p.1}
For any two vectors $v, w \in {\mathbb R}^{2d}$
\[
[Q_v,Q_w] = \frac{i}{2\pi}\Omega (v, w) {\rm Id},
\]
where the commutator $[A,B] = AB - BA$.
\end{Proposition}

For the purpose of the next definitions we shall make the following assumption: 
for given vectors $v_1, \ldots, v_d \in {\mathbb R}^{2d}$, define 
$B_{v}(j,k) = {v}^{k+d}_j$, $j,k = 1, \ldots, d$, to be a 
$d \times d$ matrix, and assume that it is {\it non-degenerate}, i.e., 
\[
\det B_{v} \neq 0.
\]

As a consequence of Proposition \ref{p.1} we observe that the
$Q_{v_j}$s commute with each other if the $v_j$s form a basis for $\Pi$.
This commutativity implies, in particular, that $\nabla_k f_j = \nabla_j f_k$, 
and so we deduce that $f_j (x) = \nabla_j (x F_{v} x)$, for some 
quadratic form $F_v$.
Thus the common eigenfunction for all the operators $\{Q_{{v}_j}\}$ 
has the form:
\[
\psi_{\xi}(x) = \frac{1}{\sqrt{|\det B_{v}}|}
e^{-2\pi i x^t B_v^{-1}\xi + 2\pi i x^tF_v x}, 
\]
for any $\xi \in {\mathbb R}^d$. 
Moreover, let $A_v(j,k) = v^k_j$, $j,k = 1, \ldots, d$. Then, the commutativity 
of the $Q_{v_j}$s implies that
\begin{equation}
\label{e4}
A_vB_v^t - B_vA_v^t = 0,
\end{equation}
where $A^t$ is the adjoint of $A$.
It follows from (\ref{e4}) that $B_v^{-1}A_v$ is symmetric. 
It is also easy to see that
\[
F_v = \frac{1}{2} B_v^{-1}A_v.
\]

Let us now define 
the following {\it gene\-ralized Fourier trans\-forms} ${\cal F}_v$
on the space of tempered distributions on ${\mathbb R}^d$,
through their action on the space of Schwartz functions: 
\[
{\cal F}_v (h) (\xi) = \int_{{\mathbb R}^d} h(x) \overline{\psi_\xi(x)} dx=
\int_{{\mathbb R}^d} h(x) \frac{1}{\sqrt{|\det B_v}|}
e^{2\pi i x^t B_v^{-1}\xi - \pi i x^tB_v^{-1}A_v x}\; dx.
\]
The operators ${\cal F}_v$ are unitary when restricted to $L^2({\mathbb R}^d)$, 
since they are combinations of unitary transformations.

We shall now consider two different representations of functions or even 
distributions 
associated with two different La\-gran\-gian planes: $\Pi$ with the basis
$v_1, \ldots, v_d$, and $\Gamma$ with the basis
$w_1, \ldots, w_d$. Assume that $\Pi \cap \Gamma =\{0\}$. Moreover, assume that
$B_v$ and $B_w$ are non-degenerate.
Define the $d\times d$ matrix $Y_{v,w} (i,j) = \Omega(v_i, w_j)$.
Note that $Y_{v,w} = {\rm Id}$ if and only if 
$\{v_1,\ldots,v_d, w_1,\ldots, w_d\}$ forms a symplectic basis for 
${\mathbb R}^{2d}$.

\begin{Lemma}
\label{t8}
\[
\det Y_{v,w} \neq 0.
\]
\end{Lemma}

\begin{Proof} Indeed, if $\Pi \cap \Gamma = \{0\}$, then
$\{v_1, \ldots, v_d, w_1, \ldots, w_d \}$ forms a 
(not necessarily symplectic) basis for ${\mathbb R}^{2d}$.
In this basis the matrix of $\Omega$ has the form:
\[
\left(\begin{array}{ll} 
0 & Y_{v,w}\\
-Y_{v,w} & 0
\end{array}\right)
\]
and so, $(\det Y_{v,w})^2 = \det \Omega \neq 0$.
\end{Proof}

The matrix $Y_{v,w}$ can be represented, with the use of 
matrices $A_v$, $A_w$, $B_v$, $B_w$, as:
\[
Y_{v,w} = A_v B_w^t - B_v A_w^t.
\]
As a consequence, we derive the following formula, which we shall use in the
proof of Theorem \ref{t12}:
\begin{equation}
\label{e5}
F_v - F_w = \frac{1}{2} \left(B_v^{-1}A_v - B_w^{-1}A_w\right) = 
\frac{1}{2} B_v^{-1} Y_{v,w} (B_w^{-1})^t. 
\end{equation}

\begin{Lemma}
\label{t9}
For any tempered distribution $h$ on ${\mathbb R}^d$, the relationship between 
its ``$v$" and ``$w$" generalized Fourier transform representations is

\begin{eqnarray*}
{\cal F}_w (h) (\eta) &=& \frac{1}{\sqrt{|\det Y_{v,w}|}} 
e^{-\pi i \eta^t Y_{v,w}^{-1} B_v B_w^{-1} \eta + \pi i \sigma/4}\\
&\times& \int_{{\mathbb R}^d} e^{\pi i \xi^t (Y_{v,w}^{-1})^t \eta 
+ \pi i \eta^t Y_{v,w}^{-1}\xi
-\pi i \xi^t (B_v^{-1})^t B_w^t Y_{v,w}^{-1} \xi}
{\cal F}_v (h)(\xi) \; d\xi,
\end{eqnarray*}
where $\sigma$ is the difference between the positive and 
negative squares of the quadratic form $F_v - F_w$.
\end{Lemma}

\begin{Proof} 
The expression in Lemma \ref{t9} is to be understood 
in the sense of distributions, 
and thus it is enough to check its validity on Schwartz functions.
Note that the inverse of the generalized Fourier transform 
${\cal F}_v$ has the form:
\[
h(x) = \frac{1}{\sqrt{|\det B_v}|} \int_{{\mathbb R}^d}
e^{-2\pi i x^t B_v^{-1} \xi + 2 \pi i x^t F_v x} {\cal F}_v (h)(\xi) \; d\xi.
\]
Taking the generalized Fourier transform ${\cal F}_w$ of this expression 
and using (\ref{e5}), we obtain
\begin{eqnarray*}
{\cal F}_w (h) (\eta)\!\!\!\! &=& \frac{1}{\sqrt{|\det B_v B_w|}} 
\int_{{\mathbb R}^d} \!
\int_{{\mathbb R}^d}\!\!\!\! e^{2\pi i x^t (B_w^{-1}\eta - B_v^{-1} \xi) + 
2 \pi i x^t (F_v - F_w) x} {\cal F}_v(h)(\xi) \; dx d\xi \\
&=&\frac{ e^{\pi i \sigma/4}}{\sqrt{|\det Y_{v,w}|}}\int_{{\mathbb R}^d} 
e^{-\pi i (B_w^{-1}\eta - B_v^{-1} \xi)^t B_w^t Y_{v,w}^{-1} B_v 
(B_w^{-1}\eta - B_v^{-1} \xi)}{\cal F}_v(h)(\xi) \;d\xi \\
&=& \frac{1}{\sqrt{|\det Y_{v,w}|}} \, e^{\pi i \sigma/4 -\pi i \eta^t 
Y_{v,w}^{-1} B_v B_w^{-1} \eta}\\
&\times& \int_{{\mathbb R}^d} e^{\pi i \xi^t (B_v^{-1})^t B_w^t Y_{v,w}^{-1} 
B_v B_w^{-1} \eta + \pi i \eta^t Y_{v,w}^{-1} \xi - \pi i \xi^t 
(B_v^{-1})^tB_w^t Y_{v,w}^{-1} \xi} {\cal F}_v(h)(\xi)\; d\xi.
\end{eqnarray*}
In order to finish the proof, it is now enough to observe that
\[
B_w^t Y_{v,w}^{-1} B_v = B_v^t (Y_{v,w}^{-1})^t B_w,
\]
due to (\ref{e5}), and that the above representation of ${\cal F}_w$ 
simplifies exactly to the formula in the statement of Lemma \ref{t9}.
\end{Proof}

We shall now introduce two more representations of 
tempered distributions associated with a collection of vectors 
$\{v_1, \ldots, v_d, w_1,\ldots, w_d\}$:
\[
{\tilde {\cal F}_v} (g)(\xi) =  e^{- \pi i \xi^t (B_v^{-1})^t
B_w^t Y_{v,w}^{-1} \xi} {\cal F}_v(g)(\xi) 
\]
and
\[
{\tilde {\cal F}_w} (g)(\eta) = e^{-\pi i \sigma/4}
e^{-\pi i \eta^t Y_{v,w}^{-1} B_v B_w^{-1} \eta}
{\cal F}_w (g)(\eta). 
\]

\noindent {\bf Remark.} In view of Steger's observation, these modifications 
of the generalized Fourier transforms
may be compared to the metaplectic representations of symplectic 
transformations which send bases of Lagrangian planes into elements of 
the standard basis for ${\mathbb R}^{2d}$.

\begin{Proposition}
\label{t10}
If $\{v_1, \ldots, v_d, w_1, \ldots, w_d \}$ forms a 
symplectic basis in ${\mathbb R}^{2d}$, i.e., $Y_{v,w} = {\rm Id}$,
then the relation between ${\tilde {\cal F}_v}$ and
${\tilde {\cal F}_w}$ takes the form of the standard Fourier transform:
\begin{equation}
\label{e8}
{\tilde {\cal F}_w}(g)(\eta) = \int_{{\mathbb R}^d}e^{2\pi i \xi \cdot \eta}
{\tilde {\cal F}_v} (g)(\xi) \; d\xi. 
\end{equation}
\end{Proposition}

\begin{Proof}
It follows easily from Lemma \ref{t9} that 
\begin{equation}
\label{e.5*}
{\tilde {\cal F}_w} (g)(\eta) =\frac{1}{\sqrt{|\det Y_{v,w}|}}
\int_{{\mathbb R}^d} e^{ \pi i [\xi^t (Y_{v,w}^{-1})^t \eta 
+ \eta^t Y_{v,w}^{-1} \xi]} 
{\tilde {\cal F}_v} (g)(\xi) \; d\xi.
\end{equation}
Since $\{v_1, \ldots, v_d, w_1, \ldots, w_d\}$ is a symplectic basis for
${\mathbb R}^{2d}$, we have $Y_{v,w} = {\rm Id}$, and so (\ref{e.5*}) 
reduces to (\ref{e8}).
\end{Proof}

We can view (\ref{e8}) as a formal and general expression for the usual 
Fourier transform of distributions. We also note that ${\tilde {\cal F}_v}$ and
${\tilde {\cal F}_w}$ are unitary transformations when 
restricted to $L^2({\mathbb R}^d)$.

It is evident that the operators ${\cal F}_v$ and ${\cal F}_w$
composed with operators $Q_{v_j}$ and $Q_{w_j}$, respectively, become 
multiplications by $j$th coordinates.
We use this fact to deduce the following lemma, which we shall use 
in the proof of our Theorem \ref{t12}.

\begin{Lemma}
\label{t11}
For each $j=1, \ldots, d$, the operators $Q_{v_j}$ are 
multiplications by $\xi_j$ in the 
${\tilde {\cal F}_v}$ representation, and all operators $Q_{w_j}$
are multiplications by $\eta_j$ in the 
${\tilde {\cal F}_w}$ representation, i.e.,
\begin{eqnarray*}
{\tilde {\cal F}_v} (Q_{v_j}(g))(\xi) &=& \xi_j {\tilde {\cal F}_v}(g) (\xi), \\
{\tilde {\cal F}_w} (Q_{w_j}(g))(\eta) &=& \eta_j {\tilde {\cal F}_w}(g) (\eta).
\end{eqnarray*}
\end{Lemma}


We can now formulate and prove our main results.

\begin{Th.}
\label{t12}
Let $\Lambda\subset {\mathbb R}^{2d}$ 
be a countable sequence of points with the property $\Lambda = - \Lambda$. 
Let $\{T_{m,n}: (m,n) \in \Lambda\}$ be the family of associated 
translation-modulation transformations $T_{m,n}$.
For a function $g \in L^2({\mathbb R}^d)$, assume that 
$\{ g_{m,n} = T_{m,n} (g): (m,n) \in \Lambda \}$ forms 
an orthonormal basis for $L^2({\mathbb R}^d)$.
For any two vectors $v,\, w \in {\mathbb R}^{2d}$ for which the symplectic form
is non-vanishing, i.e.,
\[
\Omega(v, w) \neq 0,
\]
we have
\begin{equation}
\label{e.10*}
\| Q_v (g) \|_2 \| Q_w (g) \|_2 = \infty.
\end{equation}
\end{Th.}

\begin{Proof} $i.$ Without loss of generality we may assume that 
$\Omega(v, w) =1$. There exists a collection of vectors
$\{ v_2, \ldots, v_d, w_2, \ldots, w_d \}\subset {\mathbb R}^{2d}$
such that if we let $v_1 = v$ and $w_1=w$, then 
$\{ v_1, \ldots, v_d, w_1, \ldots, w_d \}$ forms a symplectic basis
of ${\mathbb R}^{2d}$. (This result is a simple algebraic fact; for 
its Hermitian version see \cite{H}.)
With these vectors we associate the corresponding differential operators
$Q_{v_1}, \ldots, Q_{v_d}$, $Q_{w_1}, \ldots, Q_{w_d}$, and the induced 
$d\times d$ matrices $A_v, A_w, B_v, B_w$.
For this part of the proof assume that 
\[
\det B_v \neq 0 \quad {\rm and}\quad \det B_w \neq 0.
\]

Due to the assumption about the basis $\{v_1, \ldots, v_d, w_1, 
\ldots, w_d \}$, the matrix
\[
\left(\begin{array}{ll} A_v^t & B_v^t \\ 
A_w^t & B_w^t \end{array}\right)
\]
is symplectic, i.e.,
\begin{equation}
\label{e.6*}
A_v B_v^t - B_vA_v^t =0, \quad
A_w B_w^t - B_wA_w^t =0,
\end{equation}
and
\begin{equation}
\label{e.7*}
A_v B_w^t - B_vA_w^t ={\rm Id}, \quad
A_w B_v^t - B_wA_v^t =-{\rm Id}.
\end{equation}

Given a vector $(p, q) \in {\mathbb R}^{2d}$, we use translation by $x$ and
the symmetry of $B_v^{-1}A_v$ to calculate
\begin{eqnarray*}
&&{\tilde {\cal F}_v} (T_{p,q} (g)) (\xi) \\
&=& \frac{1}{|\sqrt{\det B_v}|}
e^{-\pi i \xi^t (B_v^{-1})^tB_w^t\xi} \int_{{\mathbb R}^d} 
g(x+p) e^{2\pi i x^t (B_v^{-1}\xi + q) - \pi i x^tB_v^{-1}A_v x}\; dx\\
&=& c_{p,q} e^{-\pi i\xi^t (B_v^{-1})^t B_w^t \xi -2 \pi i p^tB_v^{-1} \xi}
{\cal F}_v (g) (\xi+B_v q+ A_v p),
\end{eqnarray*}
where $c_{p,q}$ is a complex constant of absolute value equal to 1.
Recall that for a symplectic basis, $Y_{v,w} = {\rm Id}$.
Because of this and the symmetry of
$B_w B_v^{-1}$, which, in turn, follows from (\ref{e5}), we obtain
\[
{\tilde {\cal F}_v} (T_{p,q} (g)) (\xi) =
 c_{p,q} e^{2\pi i(-(B_v^{-1})^tp+B_wq +B_wB_v^{-1}A_v p)\cdot \xi}
{\tilde{\cal F}_v} (g) (\xi+B_v q+ A_v p).
\]
Therefore we can write
\[
{\tilde {\cal F}_v} (T_{p,q} (g)) (\xi) = 
c_{p,q} T_{(p', q')} \left({\tilde {\cal F}_v} (g) \right)(\xi),
\]
where 
\begin{equation}
\label{e.8*}
\left(\begin{array}{l} p' \\ q' \end{array}\right) = 
\left(\begin{array}{l} A_v p+B_v q\\
-(B_v^{-1})^tp+B_wq+B_wB_v^{-1}A_vp\end{array}\right).
\end{equation}
(\ref{e5}) yields $A_w = -(B_v^{-1})^t+B_w(B_v^{-1}A_v)^t$. 
Thus, using the symmetry of $B_v^{-1}A_v$, we can write (\ref{e.8*})
in a more familar form
\[
\left(\begin{array}{l} p' \\ q' \end{array}\right) = 
\left(\begin{array}{ll} A_v & B_v \\ 
A_w & B_w 
\end{array}\right)
\left(\begin{array}{l} p \\ q \end{array}\right).
\]
Overall, we obtain that in the ${\tilde {\cal F}}_v$ representation, 
a Gabor system remains a Gabor system, but associated with a new set $\Lambda'$:
\[
{\tilde {\cal F}_v} (T_{m,n} (g)) = c_{m',n'} ({\tilde {\cal F}_v} (g))_{m',n'},
\] 
where the primes indicate the elements of the new sequence.
Since we know that ${\tilde {\cal F}_v}$ is unitary on $L^2({\mathbb R}^d)$,
if $\{g_{m,n}: (m,n) \in \Lambda\}$ is an orthonormal basis for 
$L^2({\mathbb R}^d)$ then so is $\{c_{m',n'}({\tilde {\cal F}_v} (g))_{m',n'} : 
(m',n')\in \Lambda'\}$, where $\Lambda ' = - \Lambda'$.
Thus, using Theorem \ref{t6} and invoking Proposition \ref{t10}, we obtain that
\begin{equation}
\label{e.9*}
\| \xi_1 {\tilde {\cal F}_v} (g)(\xi) \|_2 
\| \eta_1 {\tilde {\cal F}_w} (g)(\eta) \|_2 = \infty.
\end{equation}
Moreover, because of Lemma \ref{t11}, 
we know that $Q_{v_1}$ becomes multiplication by $v_1$ in 
the ${\tilde {\cal F}_v}$ representation, and similarly  $Q_{w_1}$ becomes 
multiplication by $w_1$ in the ${\tilde {\cal F}_w}$ representation, both in
the sense of distributions.
Thus (\ref{e.9*}) is equivalent to (\ref{e.10*}), since the generalized Fourier 
transforms are unitary on $L^2({\mathbb R}^d)$ and because of 
Lemma \ref{t11}.

$ii.$ First, let us observe that, since $\Omega(v, w) =1$, we cannot have both
$(v^{1+d}, \ldots, v^{2d}) = 0$ and $(w^{1+d}, \ldots, w^{2d}) = 0$.
Therefore, without loss of generality, we assume that 
$(v^{1+d}, \ldots, v^{2d}) \neq 0$.

Recall that according to (\ref{eQ}) 
we write $Q_v = \frac{i}{2\pi} \nabla_v + f_v$.
Thus we may always find a (non-unique) non-degenerate linear transformation 
of ${\mathbb R}^d$ such that the operator $\frac{i}{2\pi} \nabla_v$
becomes $\frac{i}{2\pi} \frac{\partial}{\partial {\tilde x}^{1}}$,
in the new coordinates. The operator $Q_v$ can be
then written as
\[
Q_v = \frac{i}{2\pi} \frac{\partial}{\partial {\tilde x}^{1}} +
{\tilde a}^{1}_{v} {\tilde x}^{1} + \ldots + 
{\tilde a}^{d}_{v} {\tilde x}^{d}.
\]
We also note that
\[
{\tilde a}^{1}_{v} {\tilde x}^{1} + \dots +
{\tilde a}^{d}_{v} {\tilde x}^{d} = 
{\partial \over \partial {\tilde x}^{1}} \left(
{\tilde a}^{1}_{v} {({\tilde x}^{1})^{2} \over 2} +
{\tilde a}^{2}_{v} {\tilde x}^{1} {\tilde x}^{2} + \dots +
{\tilde a}^{d}_{v} {\tilde x}^{1} {\tilde x}^{d} \right)
= {\partial \over \partial {\tilde x}^{1}} q({\tilde x}).
\]
We define a unitary transformation $U$ of $L^2({\mathbb R}^d)$ to be
\[
U(g) (\tilde x) = e^{-2\pi i q(\tilde x)}g(\tilde x).
\]
It is easy to verify that the operator $Q_v$ takes the form 
$\frac{i}{2\pi} \frac{\partial}{\partial {\tilde x}^{1}}$ in this 
new representation, i.e.,
\[
U(Q_v(g))(\tilde x) = \frac{i}{2\pi} \frac{\partial}{\partial 
{\tilde x}^{1}} U(g) (\tilde x).
\]
Also, the operator $U \circ Q_w \circ U^{-1}$ may be written in an analogous 
form:
\[
{\tilde b}^{1}_{w} \frac{i}{2\pi}{\partial \over \partial {\tilde x}^{1}} 
+ \ldots + {\tilde b}^{d}_{w} \frac{i}{2\pi}
{\partial \over \partial {\tilde x}^{d}} + {\tilde a}^{1}_{w} {\tilde x}^{1} 
+ \ldots + {\tilde a}^{d}_{w} {\tilde x}^{d}.
\]
We shall consider three different possibilities for the differential part of 
the operator $U \circ Q_w \circ U^{-1}$.

$ii.a.$ In case $\tilde b_w = ({\tilde b}^{1}_{w},\ldots,{\tilde b}^{d}_{w}) =0$, 
$U \circ Q_w \circ U^{-1}$ has the form
\[
{\tilde a}^{1}_{w} {\tilde x}^{1} 
+ \ldots + {\tilde a}^{d}_{w} {\tilde x}^{d}.
\]
Since $\Omega (v,w) =1$, we have ${\tilde a}^{1}_{w} =1$.
We make the following non-degenerate linear transformation in ${\mathbb R}^d$:
\begin{equation}
\label{lintr2}
z^{1} = {\tilde a}^{1}_{w} {\tilde x}^{1} + \dots +
{\tilde a}^{d}_{w} {\tilde x}^{d}, \quad
z^{2} = {\tilde x}^{2} , \quad \ldots, \quad z^{d} = {\tilde x}^{d} .
\end{equation}
Thus we obtain
\[
U \circ Q_v \circ U^{-1} = \frac{i}{2\pi} \frac{\partial}{\partial z^1}, \quad
U \circ Q_w \circ U^{-1} = z^1,
\]
and the problem reduces to the standard Balian-Low theorem, Theorem \ref{t1}.

$ii.b.$ If $\tilde b_w = \alpha b_v = (\alpha, 0,\ldots, 0)$ and $\alpha 
\neq 0$ then, since we again have
${\tilde a}^{1}_{w} = 1$, by making the same transformation (\ref{lintr2})
as in part $ii.a$, we obtain:
\[
U \circ Q_v \circ U^{-1} = \frac{i}{2\pi} {\partial \over \partial z^{1}},
\quad U \circ Q_w \circ U^{-1} =
\alpha \frac{i}{2\pi} {\partial \over \partial z^{1}} + z^{1}.
\]
It is again easy to see that our result follows from the standard Balian-Low
theorem.

$ii.c.$ Finally we consider the case $\tilde b_w \neq \alpha b_v$ for all 
$\alpha$. We can make a linear transformation in 
${\mathbb R}^d$ such that $U \circ Q_v \circ U^{-1}$ remains 
the differentiation with respect to the first coordinate $z^{1}$,
and the differential part of $U \circ Q_w \circ U^{-1}$
becomes the differentiation with respect to the second coordinate $z^{2}$, i.e.,
\[
U \circ Q_v \circ U^{-1} = \frac{i}{2\pi} {\partial \over \partial z^{1}},\quad
U \circ Q_w \circ U^{-1} = \frac{i}{2\pi} {\partial \over \partial z^{2}} 
+ z^{1} + c^{2}_{w} z^{2} + \ldots + c^{d}_{w} z^{d}.
\]
We now define the following two families of vectors in ${\mathbb R}^{2d}$:
\begin{eqnarray*}
v_1 &=& (0,\ldots, 0;1, 0,\ldots, 0) \\
v_2 &=& (0, c_w^2, c_w^3,\ldots,c_w^d ;0,1,0,\ldots, 0) \\
v_3 &=& (0, c_w^3, 0, \ldots, 0;0,0,1,0,\ldots, 0) \\
&\ldots& \\
v_d &=& (0, c_w^d, 0, \ldots, 0;0,\ldots, 0, 1)
\end{eqnarray*}
and
\begin{eqnarray*}
w_1 &=& (1, c_w^2, c_w^3,\ldots,c_w^d;0,1, 0,\ldots, 0) \\
w_2 &=& (0,1,0,\ldots,0 ;1,0,\ldots, 0) \\
w_3 &=& (0, c_w^3, 1,0, \ldots, 0;0,0,1,0,\ldots, 0) \\
&\ldots& \\
w_d &=& (0,c_w^d,0,\ldots, 0, 1;0,\ldots, 0, 1),
\end{eqnarray*}
and associated with them operators
\begin{eqnarray*}
Q_{v_1} &=&  \frac{i}{2\pi} \frac{\partial}{\partial z^1}\\
Q_{v_2} &=& \frac{i}{2\pi} \frac{\partial}{\partial z^2} + c_w^2 z^2 + 
c_w^3 z^3 +\ldots +c_w^d z^d  \\
Q_{v_3} &=& \frac{i}{2\pi} \frac{\partial}{\partial z^3} + c_w^3 z^2 \\
&\ldots& \\
Q_{v_d} &=& \frac{i}{2\pi} \frac{\partial}{\partial z^3} + c_w^d z^2
\end{eqnarray*}
and
\begin{eqnarray*}
Q_{w_1} &=& \frac{i}{2\pi} \frac{\partial}{\partial z^2} + z^1+ 
c_w^2 z^2 + c_w^3 z^3 +\ldots + c_w^d z^d \\
Q_{w_2} &=& \frac{i}{2\pi} \frac{\partial}{\partial z^1} + z^2\\
Q_{w_3} &=& \frac{i}{2\pi} \frac{\partial}{\partial z^3} + c_w^3 z^2 + z^3 \\
&\ldots& \\
Q_{w_d} &=& \frac{i}{2\pi} \frac{\partial}{\partial z^d} +c_w^d z^2 + z^d .
\end{eqnarray*}
It is not difficult to verify that $\{ v_1, \ldots, v_d; w_1, \ldots, w_d\}$
forms a symplectic basis in ${\mathbb R}^{2d}$ and that the matrices 
$B_v$ and $B_w$ are both non-degenerate. Thus we have reduced this situation to
the case described in part $i$. 
\end{Proof}

\vskip0.3cm
\noindent {\bf Remark.} 
We used the notion of a symplectic matrix in the proof of Theorem \ref{t12}. 
A matrix $M$ is {\it symplectic} if it preserves 
the symplectic form $\Omega$, i.e., $\Omega(Mv, Mw) = \Omega(v,w)$,
for all $v,w \in {\mathbb R}^{2d}$. The collection of all such matrices 
forms a group, the so-called {\it symplectic group}, which plays 
a significant role in the study of 
Hamiltonian systems. In fact, the symplectic matrices generate invertible 
transformations which take a Hamiltonian system into another such system 
of differential equations, see, e.g., \cite{A}, \cite{AN}, \cite{SM}.

Following \cite{GHHK} we say that a lattice $\Lambda \subset {\mathbb R}^{2d}$
is {\it symplectic} if
\[
\Lambda = r M({\mathbb Z}^{2d})
\]
for some $r \in {\mathbb R}\setminus \{0\}$ and $M$ a symplectic
matrix. A generalized Fourier transform $\tilde {\cal F}_v$ maps a symplectic
lattice $\Lambda$ into another symplectic lattice $\Lambda'$, 
according to the formula (\ref{e.8*}).
\vskip0.3cm

\begin{Th.}
\label{t13}
Let $\Lambda\subset {\mathbb R}^{2d}$ be a lattice. 
Let $\{T_{m,n}: (m,n) \in \Lambda\}$ be the family of associated 
translation-modulation transformations $T_{m,n}$.
For a function $g \in L^2({\mathbb R}^d)$, assume that 
$\{ g_{m,n} = T_{m,n} (g): (m,n) \in \Lambda \}$ forms 
an exact frame for $L^2({\mathbb R}^d)$ and let $\tilde g$ be the canonical dual 
to $g$.
For any two vectors $v,\, w \in {\mathbb R}^{2d}$ for which the symplectic form
is non-vanishing, i.e.,
\[
\Omega(v, w) \neq 0,
\]
we have
\begin{equation}
\label{e.11*}
\| Q_v (g) \|_2 \| Q_w (g) \|_2 \| Q_v (\tilde g) \|_2 \| Q_w (\tilde g) \|_2 
= \infty.
\end{equation}
\end{Th.}

\begin{Proof} The proof is analogous to the proof of Theorem \ref{t12}.
We start with the case where vectors $v = v_1, w= w_1$ allow an extension 
$\{ v_1, \ldots, v_d, w_1, \ldots, w_d \}$ which forms a symplectic basis
of ${\mathbb R}^{2d}$ and has non-degenerate associated matrices 
$B_v$ and $B_w$. The generalized Fourier transforms $\tilde {\cal F}_v$
and $\tilde {\cal F}_w$ change the operators $Q_v$ and $Q_w$ into position 
and momentum operators, $P_1$ and $M_1$, in appropriate representations, 
respectively. Moreover, $\tilde {\cal F}_v$ maps the lattice $\Lambda$ into
another lattice $\Lambda'\subset {\mathbb R}^{2d}$. Since generalized Fourier 
transforms are unitary in $L^2({\mathbb R}^d)$, we finish by using, instead of 
Theorem \ref{t6}, a version of Theorem \ref{t1} for general lattices, see 
the remark after the proof of Theorem \ref{t1}.

The general case is reduced to the above situation analogously to the general 
case in Theorem \ref{t12}.
\end{Proof}

\vskip0.3cm

\noindent{\bf Acknowledgments.} The second named author presented this work 
at the Analysis and Probability Related to Solvable Lie Groups conference
in Zakopane in June 2002. Afterwards it was brought to our attention 
by Professor Tim Steger that Theorem \ref{t12} can also be proved using 
the theory of metaplectic representations.
We gratefully acknowledge Professor Steger for his comments.
We also gratefully acknowledge Professor Raymond Johnson for making us aware 
of the connection between the work of H\"ormander and our results.

The first named author was partially supported by NSF-DMS Grant 0139759 
(2002--2005) and the General Research Board of the University of Maryland.


\newpage
\noindent \author{John J. Benedetto\\
                Department of Mathematics\\
                University of Maryland\\
                College Park, MD 20742, USA\\
                {\tt jjb@math.umd.edu}}
\vskip1cm                
\noindent \author{Wojciech Czaja\\
                Mathematical Institute\\
                University of Wroc\l aw\\
                50-384 Wroc\l aw, pl. Grunwaldzki 2/4, Poland\\
                {\tt czaja@math.uni.wroc.pl}\\
                and\\
                Department of Mathematics\\
                University of Maryland\\
                College Park, MD 20742, USA\\
                {\tt wojtek@math.umd.edu}}
\vskip1cm       
\noindent \author{Andrei Ya. Maltsev\\
                L.D.Landau Institute for Theoretical Physics\\
                119334 ul. Kosygina 2, Moscow, Russia\\
                {\tt maltsev@itp.ac.ru}\\
                and\\
                Department of Mathematics\\
                University of Maryland\\
                College Park, MD 20742, USA}

\end{document}